\documentclass[11pt]{article}
\usepackage{graphicx}
\usepackage{amsmath,amsfonts,amssymb}
\usepackage{citesort}
\usepackage{url}
\usepackage{mdframed}
\usepackage{amsmath}
\usepackage{graphicx}
\usepackage{epsfig}
\usepackage{epstopdf}
\usepackage{color}
\def\Limsup{\mathop{{\rm Lim}\,{\rm sup}}}

\def\tto{\;{\lower 1pt \hbox{$\rightarrow$}}\kern -10pt
\hbox{\raise 2pt \hbox{$\rightarrow$}}\;}

\def\argmin{\mathop{\rm argmin}\nolimits}
\def\ra{\rangle}
\def\la{\langle}

\def\B{\Bbb B}

\def\qed{$\hfill\Box$}
\def\R{\Bbb R}
\def\N{\Bbb N}
\def\ox{\bar{x}}
\def\ow{\bar{w}}

\def\oz{\bar{z}}

\def\gph{\mbox{\rm gph}\,}

\def\dom{\mbox{\rm dom}\,}

\def\O{\Omega}

\newcounter{lk}

\begin{document}

\begin{center}
\vspace*{0.3in} {\bf CONVERGENCE ANALYSIS OF A PROXIMAL POINT ALGORITHM FOR MINIMIZING DIFFERENCES OF FUNCTIONS}
\\[2ex]
 Nguyen Thai An\footnote{Thua Thien Hue College of Education, 123 Nguyen Hue, Hue City, Vietnam (thaian2784@gmail.com). The research of Nguyen Thai An was supported by the Vietnam National Foundation for Science and Technology Development (NAFOSTED) under grant number 101.01-2014.37.},
Nguyen Mau Nam\footnote{Fariborz Maseeh Department of Mathematics and Statistics, Portland State University, PO Box 751, Portland, OR 97207, United States (mau.nam.nguyen@pdx.edu). The research of Nguyen Mau Nam was partially supported by the USA National Science Foundation under grant DMS-1411817 and
the Simons Foundation under grant \#208785.},
\\[1ex]
\end{center}
{\small \textbf{Abstract.} Several optimization schemes have been known for convex optimization problems. However, numerical algorithms for solving nonconvex optimization problems are still underdeveloped. A progress to go beyond convexity was made by considering the class of functions representable as differences of convex functions. In this paper, we introduce a generalized proximal point algorithm to minimize the difference of a nonconvex function and a convex function. We also study convergence results of this algorithm under the main assumption that the objective function satisfies the Kurdyka - \L ojasiewicz property.

%The main goal of this paper is to analyze the convergence of a new proximal point algorithm (GPPA) introduced by Souza et al.  for minimizing differences of convex functions. }

\medskip
\vspace*{0,05in} \noindent {\bf Key words.} DC programming,  proximal point algorithm, difference of convex functions, Kurdyka - \L ojasiewicz
inequality.

\newtheorem{Theorem}{Theorem}[section]
\newtheorem{Proposition}[Theorem]{Proposition}
\newtheorem{Remark}[Theorem]{Remark}
\newtheorem{Lemma}[Theorem]{Lemma}
\newtheorem{Corollary}[Theorem]{Corollary}
\newtheorem{Definition}[Theorem]{Definition}
\newtheorem{Example}[Theorem]{Example}
\renewcommand{\theequation}{\thesection.\arabic{equation}}
\normalsize

\section{Introduction}
In this paper, we introduce and study the convergence analysis of an algorithm for solving optimization problems in which the objective functions can be represented as differences of nonconvex and convex functions. The structure of the problem under consideration is flexible enough to include the problem of minimizing a smooth function on a closed set or minimizing a DC function, where DC stands for {\bf D}ifference of {\bf C}onvex functions. It is worth noting that DC programming is one of the most successful approaches to go beyond convexity. The class of DC functions is closed under many operations usually considered in optimization and is quite large to contain many objective functions in applications of optimization. Moreover, this class of functions possesses beautiful generalized differentiation properties and is favorable for applying numerical optimization schemes; see \cite{Tuy98,HorTuy90,BacBor11} and the references therein.

A pioneer in this research direction is Pham Dinh Tao who introduced a simple algorithm called the (DCA) based on generalized differentiation of the functions involved as well as their Fenchel conjugates \cite{TaoSou86}. Over the past three decades, Pham Dinh Tao, Le Thi Hoai An and many others have contributed to providing mathematical foundation for the algorithm and making it accessible for applications. The (DCA) nowadays becomes a classical tool in the field of optimization due to several key features including simplicity, inexpensiveness, flexibility and efficiency; see \cite{TA1,TA2,TaoAnAkoa05,MuuQuoc10}.

The {\em proximal point algorithm} (PPA for short) was suggested by Martinet \cite{Mar72} for solving convex optimization problems and was extensively developed by Rockafellar \cite{Rock76} in the
context of monotone variational inequalities. The main idea of this method consists of replacing the initial problem with
a sequence of regularized problems, so that each particular auxiliary problem can be solved by one of the well-known algorithms. Along with the (DCA), a number of proximal point optimization schemes have been proposed in \cite{SSC03,MM06,BFO15,SONS15} to minimize differences of convex functions. Although  convergence results for the (DCA) and the proximal point algorithms for minimizing differences of convex functions have been addressed in some recent research, it is still an open research question to study the convergence analysis of algorithms for minimizing differences of functions in which convexity is not assumed.

Based on the method developed recently in \cite{ABRS10,BE2015,TAN}, we study a proximal point algorithm for minimizing the difference of nonsmooth functions in which only the second function involved is required to be convex. Under the main assumption that the objective function satisfies the Kurdyka - \L ojiasiewicz property, we are able to analyze the convergence of the algorithm. Our results further recent progress in using the Kurdyka - \L ojiasiewicz property and variational analysis to study nonsmooth numerical algorithms pioneered by Attouch, Bolte, Redont, Soubeyran, and many others. The paper is organized as follows. In Section 2, we provide tools of variational analysis used throughout the paper. Section 3 is the main section of the paper devoted to the generalized proximal point algorithm and its convergence results. Applications to trust-region subproblems and nonconvex feasibility problems are introduced in Section 4.

\section{Tools of Variational Analysis}
\setcounter{equation}{0}

In this section, we recall some basic concepts and results of generalized differentiation for nonsmooth functions used throughout the paper; see, e.g., \cite{Cl90,m-book1,bmn,r} for more details. We use $\R^n$ to denote the $n$ - dimensional Euclidean space, $\la \cdot, \cdot \ra$ to denote the inner product, and $\|\cdot\|$ to denote the associated Euclidean norm. For an extended-real-value function $f:\R^n \to \R\cup\{+\infty\}$, the domain of $f$ is the set $$\mbox{dom} f=\{x\in \R^n:\; f(x)<+\infty\}.$$
The function $f$ is said to be \emph{proper} if its domain is nonempty.

Given a lower semicontinuous function $f: \R^n \to \R\cup\{+\infty\}$ with $\ox\in \dom f$, the \emph{Fr\'echet subdifferential} of $f$ at $\ox$ is defined by
\begin{equation*}
\partial^F f(\ox)=\left\{v\in \R^n:\; \liminf_{x\to \ox}\dfrac{f(x)-f(\ox)-\la v, x-\ox\ra}{\|x-\ox\|}\geq 0\right\}.
\end{equation*}
 We set $\partial^F f(\ox)=\emptyset$ if $\ox \notin \dom f$. Note that the Fr\'{e}chet subdifferential mapping does not have a closed graph, so it is unstable computationally. Based on the Fr\'echet subdifferential, the \emph{limiting/Mordukhovich subdifferential} of $f$ at $\ox\in\dom f$ is defined by
\begin{equation*}
\partial^Lf(\ox)=\Limsup_{x\xrightarrow{f}\ox}\partial^F f(x)=\{v\in \R^n:\; \exists\; x^k\xrightarrow{f}\ox,  v^k\in \partial^F f(x^k), v_k \to v\}.
\end{equation*}
where the notation $x\xrightarrow{f}\ox$ means that $x\to \ox$ and $f(x)\to f(\ox)$. We also set $\partial^L f(\ox)=\emptyset$ if $\ox \notin \dom f$. It follows from the definition the following robustness/closedness  property of $\partial^L f$:
$$\bigg\{v\in \R^n:\; \exists\; x^k\xrightarrow{f}\ox, v^k\to v, v^k\in \partial^L f(x^k)\bigg\} = \partial^L f(\ox).$$
Obviously, we have $\partial^F f(x) \subset \partial^L f(x)$ for every $x\in \R^n$, where the first set is closed and convex while the second one is closed; see \cite[Theorem 8.6, p 302]{RocWet98}. If $f$ is differentiable at $\ox$, then $\partial^F(\ox)=\{\nabla f(\ox)\}$. Moreover, if $f$ is continuously differentiable on a neighborhood of $\ox$, then $\partial^L f(\ox)=\{\nabla f(\ox)\}$. When $f$ is convex, the Fr\'{e}chet and the  limiting subdifferentials reduce to the subdifferential in the sense of convex analysis:
$$\partial f(\ox)=\{v\in \R^n:\; \la v, x-\ox\ra \leq f(x) -f(\ox), \; \forall \;x\in \R^n\}.$$
For a convex subset $\O$ of $\R^n$ and $\ox\in \O$, the \emph{normal cone} to $\O$ at $\ox$ is the set
\begin{equation*}
N(\ox; \O)=\{v\in \R^n:\; \la v, x-\ox\ra\leq 0, \; \forall \;x\in \O\}.
\end{equation*}
This normal cone can be represented as the subdifferential at the point under consideration of the \emph{indicator function}:
\begin{equation*}\label{sub_dis}
   \delta (x; \O)= \left\{
     \begin{array}{lr}
       0 & \mbox{\rm if }\;\ox \in \O,\\
       +\infty & \mbox{\rm if }\; \ox \notin \O,
     \end{array}
   \right.
\end{equation*}
i.e., $N(\ox; \O)=\partial \delta(\ox; \O)$. We use the notation $\mbox{dist}(\ox; \O)$ to denote the distance from $\ox$ to $\O$, i.e., $\mbox{dist}(\ox; \O)=\inf_{x\in \O}\|x-\ox\|$. The notation $P_{\O}(\ox)=\{\ow\in \O: \; \|\ox-\ow\|=\mbox{dist}(\ox; \O)\}$ stands for the projection from $\ox$ onto $\O$. We also use $d_\O(\ox)$ for $\mbox{dist}(\ox; \O)$ where convenience.

Another subdifferential concept called the \emph{Clarke subdifferential} was defined in  \cite{Cl90} based on generalized directional derivatives. The Clarke subdifferential of a locally Lipschitz continuous function  $f$ around $\ox$ can be represented in terms of the limiting subdifferential:
\begin{equation*}
\partial^Cf(\ox)={\rm co}\,\partial^Lf(\ox).
\end{equation*}
Here ${\rm co}\,\O$ denotes the convex hull of an arbitrary set $\O$.

\begin{Proposition}{\rm (\cite[Exercise 8.8, p. 304]{RocWet98}).} \label{sumrule} Let $f=g+h$ where $g$ is lower semicontinuous and let $h$ is continuously differentiable on a neighborhood of $\ox$. Then
$$\partial^F f(\ox)=\partial^F g(\ox) +\nabla h(\ox)\; \mbox{ and }\; \partial^L f(\ox)=\partial^L g(\ox) +\nabla h(\ox).$$
\end{Proposition}

\begin{Proposition}{\rm (\cite[Theorem 10.1, p. 422]{RocWet98}).} \label{optcond}
If a lower semicontinuous function $f: \R^n \to \R \cup \{+\infty\}$ has a local minimum at $\ox\in \dom f$, then
$0\in \partial^F f(\ox) \subset \partial^L f(\ox)$. In the convex case, this condition is not only necessary for a local minimum
but also sufficient for a global minimum.
\end{Proposition}

\begin{Proposition}\label{subbound}
Let $h:\R^n \to \R$ be a finite convex function on $\R^n$. If $y^k \in \partial h(x^k)$ for all $k$ and  $\{x^k\}$ is bounded, then the sequence $\{y^k\}$ is also bounded.
\end{Proposition}
{\it Proof}  The result follows from the fact that $h$ is locally Lipschitz continuous on $\R^n$ and \cite[Definition 5.14, Proposition 5.15, Theorem 9.13]{RocWet98}. \qed

%It is easy to check that \eqref{findxk} is equavalent to requirement that
% $$y^k \in \partial g(x^{k+1}) + \frac{1}{\lambda}(x^{k+1} -x^k).$$
% This is equivalent to saying that $x^{k+1} \in \partial \Tilde{g}^*(y^k)$, where $\Tilde{g}(x) = g(x) + \frac{1}{2\lambda}\|x-x^k\|^2$.\\[1ex]

Following \cite{ABRS10,BE2015}, a lower semicontinuous function $f\colon \R^n\to \R\cup\{+\infty\}$ satisfies the \emph{Kurdyka - \L ojasiewicz property} at $x^*\in \mbox{\rm dom}\,\partial^Lf$ if there exist $\nu>0$, a neighborhood $\mathcal{V}$ of $x^*$, and a continuous concave function $\varphi: [0, \nu[ \to [0, +\infty[$ with \\[1ex]
{\bf(i)} $\varphi(0)=0$.\\
{\bf(ii)} $\varphi$ is of class $C^1$ on $]0, \nu[$.\\
{\bf(iii)} $\varphi'>0$ on $]0, \nu[$.\\
{\bf(iv)} For every $x\in \mathcal{V}$ with $f(x^*)< f(x) < f(x^*)+\nu$, we have
$$\varphi'\left(f(x) - f(x^*)\right)\mbox{dist}\left(0, \partial^L f(x)\right)\geq 1.$$
We say that $f$ satisfies the \emph{strong Kurdyka - \L ojasiewicz} property at $x^*$ if the same assertion  holds for the Clarke subdifferential $\partial^Cf(x)$.

According to \cite[Lemma 2.1]{ABRS10}, a  proper lower semicontinuous function $f: \R^n\to \R\cup\{+\infty\}$ has the Kurdyka - \L ojasiewicz property at any point $\ox\in \R^n$ such that $0\notin \partial^L f(\ox)$. Recall that a subset $\O$ of $\R^n$ is called \emph{semi-algebraic} if it can be represented as a finite union of sets of the form
\begin{equation*}
\{x\in \R^n:\; p_i(x)=0, q_i(x)<0\; \mbox{\rm for all }i=1, \ldots, m\},
\end{equation*}
where $p_i$ and $q_i$ for $i=1,\ldots, m$ are polynomial functions.  A function $f$ is said to be semi-algebraic if its graph is a semi-algebraic subset of $\R^{n+1}$. It is known that a proper lower semicontinuous semi-algebraic function always satisfies the Kurdyka - \L ojasiewicz property; see \cite{ABRS10,BDLS07}. In a recent paper, Bolte et al. \cite[Theorem 14]{BDLS07} showed that the class of \emph{definable functions}, which contains the class of semi-algebraic functions, satisfies the strong Kurdyka - \L ojasiewicz property at each point of $\dom \partial^C f$.

\section{A Generalized Proximal Point Algorithm for Minimizing Differences of Functions}

We focus on the convergence analysis of a proximal point algorithm for solving nonconvex optimization problems of the following type
\begin{equation}
\min \bigg\{f(x) = g_1(x) + g_2(x) - h(x):\; x\in \R^n\bigg\},
\label{generalDC}\end{equation}
where $g_1(x)\colon \R^n\to \R\cup\{+\infty\}$ is proper and lower semicontinuous, $g_2(x)\colon \R^n \to \R$ is differentiable with $L$ - Lipschitz gradient, and $h: \R^n \to \R$ is convex. The specific structure of \eqref{generalDC} is flexible enough to include the problem of minimizing a smooth function on a closed constraint set:
\begin{equation*}
\min\{g(x):\; x\in \O\},
\end{equation*}
and the general DC problem:
\begin{equation}
\min\big\{ f(x) =g(x) -h(x): x\in \R^n\big\},
\label{d.c.}
\end{equation}
where $g\colon \R^n \to \R\cup\{+\infty\}$ is a proper lower semicontinuous convex function and $h\colon \R^n\to \R$ is convex.

It is well-known that if $\ox\in \dom f$ is a local minimizer of ~\eqref{d.c.}, then
\begin{equation}\label{dcstationary}
\partial h(\ox) \subset \partial g(\ox).
\end{equation}
Any point $\ox\in \dom f$ that satisfies \eqref{dcstationary} is called a \emph{stationary point} of ~\eqref{d.c.}, and any point $\ox\in \dom f$ such that $\partial g(\ox)\cap \partial h(\ox)\neq \emptyset$ is called a \emph{critical point} of this problem. Since $h$ is a finite convex function, its subdifferential at any point is nonempty, and hence any stationary point of  \eqref{d.c.} is a critical point; see \cite{TA1,HU85,HorstThoai99} and the references therein for more details.

Let us recall below a necessary optimality condition from \cite{bny} for minimizing differences of functions in the nonconvex settings.

\begin{Proposition}{\rm (\cite[Proposition~4.1]{bny})}\label{bny} Consider the difference function $f=g-h$, where $g: \R^n\to \R\cup\{+\infty\}$ and $h: \R^n \to \R$ are lower semicontinuous functions. If $\ox\in \dom f$ is a local minimizer of $f$, then we have the inclusion
$$\partial^Fh(\ox) \subset \partial^F g(\ox).$$
If in addition $h$ is convex, then  $\partial h(\ox) \subset \partial^L g(\ox).$
\end{Proposition}
When adapting to the setting of \eqref{generalDC}, we obtain the following optimality condition.

\begin{Proposition}\label{DCO}
If $\ox\in \mbox{\rm dom}\,f$ is a local minimizer of the function $f$ considered in~\eqref{generalDC}, then
\begin{equation} \label{stationary}
\partial h(\ox) \subset \partial^L g_1(\ox) + \nabla g_2(\ox).
\end{equation}
\end{Proposition}
{\it Proof } The assertion follows from Proposition \ref{sumrule} and Proposition \ref{bny}.\qed

Following the DC case, any point $\ox\in \dom f$ satisfying condition~\eqref{stationary} is called a stationary point of~\eqref{generalDC}. In general, this condition is hard to be reached and we may relax it to
\begin{equation}
[\partial^L g_1(\ox) +\nabla g_2(\ox)] \cap \partial h(\ox)\neq \emptyset
\label{criticalpoint}
\end{equation}
and call $\ox$ a  critical point of $f$. Obviously, every stationary point $\ox$ is a critical point. Moreover, by \cite[Corollary 3.4]{bny} at any point $\ox$ with $g_1(\ox)<+\infty$, we have
$$\partial^L(g_1 + g_2 -h)(\ox) \subset \partial^Lg_1(\ox)+ \nabla g_2(\ox) -\partial h(\ox).$$
Thus, if $0\in \partial^L f(\ox)$, then $\ox$ is a critical point of $f$ in the sense of \eqref{criticalpoint}. The converse is not true in general as shown by the following example. Consider the functions below
$$f(x)=2|x| + 3x,\; g_1(x) = 3|x|, \; g_2(x)= 3x,\; \mbox{\rm and } h(x)=|x|.$$
In this case, $\ox=0$ satisfies \eqref{criticalpoint} but $0\notin \partial^L f(\ox)$ since
$\partial g_1(0)=[-3, 3]$, $\nabla g_2(0)=3$, $\partial h(0)=[-1, 1]$ and $\partial f(0)=[1, 5]$. However, it is easy to check that these two conditions are equivalent when $h$ is differentiable on $\R^n$.

We recall now the \emph{Moreau/Moreau-Yoshida proximal mapping} for a nonconvex function; see \cite[page 20]{RocWet98}. Let $g: \R^n \to \R\cup\{+\infty\}$ be a proper lower semicontinuous function. The Moreau proximal mapping with regularization parameter $t>0$, $\mbox{\bf prox}_t^g: \; \R^n \to 2^{\R^n}$,  is defined by
$$\mbox{\bf prox}_t^g(x)=\argmin\bigg\{g(u) +\frac{t}{2}\|u-x\|^2: \; u\in \R^n\bigg\}.$$
As an  interesting case, when $g$ is the indicator function $\delta(\cdot; \O)$ associated with a nonempty closed set $\O$, $\mbox{\bf prox}_t^g(x)$ coincides with the projection mapping.

Under the assumption $\inf_{x\in \R^n} g(x) > -\infty$, the lower semicontinuity of $g$ and the coercivity of the squared norm imply that the proximal mapping is well-defined; see \cite[Proposition 2.2]{BST}.

\begin{Proposition}Let $g: \R^n \to \R\cup\{+\infty\}$ be a proper lower semicontinuous function with $\inf_{x\in \R^n}g(x)>-\infty$. Then, for every $t \in (0, +\infty)$, the set $\mbox{\bf prox}_{t}^g(x)$ is nonempty and compact for every $x\in \R^n$.
\end{Proposition}

We now introduce a new generalized proximal point algorithm for solving~\eqref{generalDC}. Let us begin with the lemma below regarding an upper bound for a smooth function with Lipschitz continuous gradient; see \cite{Nes04,OR70}.

%Throughout the paper, we assume that $g_1$ is lower semicontinuous and bounded from below, $g_2$ is differentiable with $L-$Lipschitz gradient, and $h$ is convex.

\begin{Proposition}If $g: \R^n \to \R$ is a differentiable function with $L$ - Lipschitz gradient, then
\begin{equation}\label{upb}
g(y) \leq g(x) +\la\nabla g(x), y-x\ra + \frac{L}{2}\|y-x\|^2 \mbox{ for all } x, y \in \R^n.
\end{equation}
\end{Proposition}
%We say that $x^*$ is a critical point of~\eqref{generalDC} if
%$$[\partial^L g_1(x^*) + \nabla g_2(x^*)] \cap \partial h(x^*) \neq \emptyset.$$

Let us introduce the generalized proximal point algorithm (GPPA) below to solve \eqref{generalDC}.

\begin{mdframed}

{\bf Generalized Proximal Point Algorithm (GPPA)}

1. Initialization: Choose $x^0\in \mbox{dom}\; g_1$ and a tolerance $\epsilon>0$.  Fix any $t>L$.

2. Find $$y^k\in \partial h(x^k).$$

3. Find $x^{k+1}$ as follows
\begin{equation}
x^{k+1}\in \mbox{\bf prox}_t^{g_1}\left(x^k - \frac{\nabla g_2(x^k) -y^k}{t}\right).
\label{proxxk1}
\end{equation}

4. If $\|x^{k}-x^{k+1}\| \leq \epsilon$, then exit. Otherwise, increase $k$ by 1 and go back to step 2.

\end{mdframed}

From the definition of proximal mapping,~\eqref{proxxk1} is equivalent to saying that
\begin{equation}
x^{k+1}\in \argmin \limits_{x\in \R^n}\bigg\{g_1(x) - \la  y^k - \nabla g_2(x^k), x-x^k\ra  + \frac{t}{2} \|x-x^k \|^2 \bigg\}.
\label{xk1}
\end{equation}

\begin{Theorem}\label{GPPA_main}
Consider the {\rm(GPPA)} for solving~\eqref{generalDC} in which $g_1(x)\colon \R^n\to \R\cup\{+\infty\}$ is proper and lower semicontinuous with $\inf_{x\in \R^n}g_1(x)>-\infty$, $g_2(x)\colon \R^n \to \R$ is differentiable with $L$ - Lipschitz gradient, and $h: \R^n \to \R$ is convex.  Then\\
{\bf (i)} For any $k\geq 1$, we have
\begin{equation}
f(x^k) - f(x^{k+1}) \geq \frac{t-L}{2}\|x^{k}-x^{k+1}\|^2.
\label{descentf}
\end{equation}
{\bf (ii)} If $\alpha = \inf \limits_{x\in \R^n} f(x) >-\infty$, then $\lim \limits_{k\to +\infty} f(x^k)=\ell^* \geq \alpha$ and $\lim \limits_{k\to +\infty}\|x^k - x^{k+1}\|=0$.\\
{\bf (iii)} If  $\alpha = \inf \limits_{x\in \R^n} f(x) >-\infty$ and $\{x^k\}$ is bounded, then every cluster point of $\{x^k\}$ is a critical point of $f$.
\end{Theorem}
{\it Proof} {\bf (i)} By Proposition \ref{sumrule} and Proposition \ref{optcond}, it follows from \eqref{xk1} that
\begin{equation}
y^k - \nabla g_2(x^k) \in \partial^L g_1(x^{k+1}) + t\left(x^{k+1} -x^k\right).
\label{findyk}
\end{equation}
Since $y^k\in \partial h(x^k)$,
\begin{equation}
h(x^{k+1}) \geq h(x^k) + \la y^k, x^{k+1} -x^k\ra.
\label{eq10}
\end{equation}
From~\eqref{xk1}, we have
\begin{equation}
g_1(x^k) \geq g_1(x^{k+1}) - \la y^k -\nabla g_2(x^k), x^{k+1} -x^k\ra +\frac{t}{2}\|x^{k}-x^{k+1} \|^2.
\label{eq11}
\end{equation}
Adding~\eqref{eq10} and~\eqref{eq11} and using \eqref{upb}, we get
\begin{align*}
g_1(x^k) - h(x^k)&\geq g_1(x^{k+1}) - h(x^{k+1}) +\la \nabla g_2(x^{k}), x^{k+1} -x^k\ra + \frac{t}{2}\|x^k -x^{k+1}\|^2\\
&\geq g_1(x^{k+1}) - h(x^{k+1}) + \left(g_2(x^{k+1}) - g_2(x^k) - \frac{L}{2}\|x^{k}-x^{k+1}\|^2 \right)+\frac{t}{2}\|x^{k}-x^{k+1}\|^2.
\end{align*}
This implies
$$f(x^k) - f(x^{k+1}) \geq \frac{t-L}{2}\|x^k - x^{k+1}\|^2.$$
Assertion {\bf(i)} has been proved. \\[1ex]
{\bf (ii)} It follows from the assumptions made and {\bf (i)} that $\{f(x_k)\}$ is monotone decreasing and bounded below, so the first assertion of {\bf (ii)} is obvious. Observe that
$$\sum_{k=1}^{m}\|x^k-x^{k+1}\|^2\leq \frac{2}{t-L}\left(f(x^1) -f(x^{m+1})\right) \leq \frac{2}{t-L}\left(f(x^1) - \alpha\right)\; \mbox{\rm for all }m\in \N.$$
Thus, the sequence $\{\|x^{k}-x^{k+1}\|\}$ converges to $0$. \\[1ex]
{\bf (iii)} From~\eqref{xk1}, for all $x\in \R^n$, we have
\begin{equation}
 g_1(x^{k+1}) - \la w^k, x^{k+1}-x^k\ra  + \frac{t}{2}\|x^{k+1}-x^k \|^2 \leq g_1(x) - \la  w^k, x-x^k\ra  + \frac{t}{2}\|x-x^k \|^2,
\label{limest}
\end{equation}
where $w^k=y^k - \nabla g_2(x^k).$ Now suppose further that $\{x^k\}$ is bounded. Since $h$ is finite convex function on $\R^n$, $y^k\in \partial h(x^k)$ and $\{x^k\}$ is bounded, from Proposition \ref{subbound}, $\{y^k\}$ is also bounded. We can take two subsequences: $\{x^{k_\ell}\}$ of $\{x^k\}$ and $\{y^{k_\ell}\}$ of $\{y^k\}$ that converge to $x^*$ and $y^*$, respectively. Because $\|x^{k_\ell}- x^{k_\ell+1}\| \to 0$ as $\ell\to +\infty$, we deduce from~\eqref{limest} that
$$\limsup \limits_{\ell \to +\infty} g_1(x^{k_\ell+1}) \leq g_1(x) - \la y^{*} - \nabla g_2(x^{*}), x - x^{*} \ra  + \frac{t}{2} \|x - x^{*}\|^2\; \mbox{\rm for all}\; x\in \R^n.$$
In particular, for $x=x^*$, we get
$$\limsup\limits_{\ell \to +\infty} g_1(x^{k_\ell+1}) \leq g_1(x^*).$$
Combining this with the lower semicontinuity of $g_1$, we get
$$\lim \limits_{\ell \to +\infty}g_1(x^{k_\ell+1}) = g_1(x^*).$$
 From the closed property of the subdifferential mapping $\partial h(\cdot)$, we have $y^* \in \partial h(x^*)$. It follows from~\eqref{findyk} that there exists $z^{k_\ell+1}\in \partial^L g_1(x^{k_\ell+1})$ satisfying
$$\|y^{k_\ell} - \nabla g_2(x^{k_\ell}) - z^{k_\ell +1}\| =t\|x^{k_\ell} - x^{k_\ell +1}\|.$$
By {\bf(ii)} and the Lipschitz continuity of $\nabla g_2$,
$$\lim \limits_{\ell\to +\infty}z^{k_\ell+1}=y^*-\nabla g_2(x^*):=z^*.$$
Thus, $x^{k_\ell+1}\xrightarrow{g_1}x^*$,  $z^{k_\ell+1} \in \partial^L g_1(x^{k_\ell+1})$, $z^{k_\ell+1} \to z^*$ as $\ell \to +\infty$, it follows from the robustness of limiting subdifferential that $z^*\in \partial^L g_1(x^*)$. Therefore, $$y^* \in [\partial^L g_1(x^*) +\nabla g_2(x^*)] \cap \partial h(x^*).$$ This implies that $x^*$ is a critical point of $f$ and the proof is complete. \qed

\begin{Proposition}\label{flsc}
Suppose that $\inf_{x\in \R^n}f(x)>-\infty$, $f$ is proper and lower semicontinuous. If the {\rm (GPPA)} sequence $\{x^k\}$ has a cluster point $x^*$, then $\lim \limits_{k \to +\infty} f(x^{k})= f(x^*)$. Thus, $f$ has the same value at all cluster points of $\{x^k\}$.
\end{Proposition}
{\it Proof} Since $\inf_{x\in \R^n} f(x) >-\infty$, it follows from \eqref{descentf} that the sequence of real numbers $\{f(x^k)\}$  is non-increasing and bounded below. Thus, $\lim \limits_{k\to +\infty}f(x^k)=\ell^*$ exists. If $\{x^{k_\ell}\}$ is a subsequence  converging to $x^*$, then by the lower semicontinuity of $f$, we have $\liminf \limits_{\ell \to +\infty} f(x^{k_\ell}) \geq f(x^*)$. Observe from the structure of $f$ that $\dom f=\dom g_1$. Since $g_2$ and $h$ are continuous, $f$ is proper and lower semicontinuous if and only if $g_1$ is proper and lower semicontinuous. To prove the opposite inequality, we employ the proof of {\bf(iii)} of Theorem \ref{GPPA_main} and get
\begin{align*}
\limsup \limits_{\ell \to +\infty} f(x^{k_\ell}) & = \limsup \limits_{\ell \to +\infty} \left(g_1(x^{k_\ell}) + g_2(x^{k_\ell}) - h(x^{k_\ell}) \right)\\
& \leq \limsup \limits_{\ell \to +\infty} g_1(x^{k_\ell}) + \limsup \limits_{\ell \to +\infty} g_2(x^{k_\ell}) - \liminf \limits_{\ell \to +\infty} h(x^{k_\ell}) \\
&\leq g_1(x^*) + g_2(x^*) - h(x^*) = f(x^*).
\end{align*}
Combining this with the uniqueness of limit, we have $\ell^*=f(x^*)$. The proof is complete. \qed

\begin{Remark}{\rm {\bf (i)} If $g$ is also convex, we can get a stronger inequality than
\eqref{descentf} and relax the range of the regularization parameter $t$. Indeed, using definition of the subdifferential in the sense of convex analysis in \eqref{findyk}, we have
\begin{equation*}
\la y^k - \nabla g_2(x^k) -t(x^{k+1} -x^k), x^k-x^{k+1}\ra \leq g_1(x^k) -g_1(x^{k+1}).
%\label{findyknew}
\end{equation*}
Since $y^k\in \partial h(x^k)$,
\begin{equation*}
h(x^{k+1}) \geq h(x^k) + \la y^k, x^{k+1} -x^k\ra.
%\label{eq10new}
\end{equation*}
Adding these inequalities and using \eqref{upb} give
$$f(x^k) - f(x^{k+1}) \geq     \left(t -\frac{L}{2} \right)   \|x^k - x^{k+1}\|^2.$$
Thus, we can choose $t>\frac{L}{2}$ instead of $t>L$ as before.\\[1ex]
{\bf(ii)} When $h(x)=0$, the (GPPA) reduces to the proximal forward - backward algorithm for minimizing $f= g_1 +g_2$ considered in
\cite{ABS11}. If $h(x)=0$ and $g_1$ is the indicator function $\delta(\cdot; \O)$ associated with a nonempty closed set $\O$, then the (GPPA) reduces to the
projected gradient  method (PGM) for minimizing the smooth function $g_2$ on a nonconvex constraint set $\O$:
$$x^{k+1}=P_{\O}\left(x^k - \frac{1}{t}\nabla g_2(x_k)\right).$$
{\bf (iii)} If $g_2 = 0$, then the (GPPA) reduces to the (PPA) with constant stepsize proposed in \cite{SSC03,SOS15}.
}
\end{Remark}

In the theorem below, we establish sufficient conditions that guarantee the convergence of the sequence $\{x_k\}$ generated by the (GPPA). These conditions include the Kurdyka - \L ojasiewicz property of the function $f$ and the differentiability with Lipschitz gradient of $h$. In what follows, let $C^*$ denote the set of cluster points of the sequence $\{x^k\}$. We follow the method from \cite{ABRS10,BE2015}.

\begin{Theorem}\label{main}
Suppose that $\inf_{x\in \R^n}f(x) >-\infty$, and $f$ is lower semicontinuous. Suppose further that $\nabla h$ is $L(h)$ - Lipschitz continuous and $f$ has the Kurdyka - \L ojasiewicz property at any point $x\in \mbox{\rm dom} f$.  If $C^*\neq \emptyset$, then the  {\rm (GPPA)} sequence $\{x^k\}$ converges to a critical point of $f$.
\end{Theorem}
{\it Proof }Take any $x^*\in C^*$ and a subsequence $\{x^{k_\ell}\}$ that converges to $x^*$. Applying Proposition \ref{flsc} yields
$$\lim \limits_{k\to +\infty}f(x^k)=\ell^*=f(x^*).$$
If $f(x^k)=\ell^*$ for some $k\geq 1$, then $f(x^k)=f(x^{k+p})$ for any $p\geq 0$ since the sequence $\{f(x^k)\}$ is monotone decreasing by \eqref{descentf}. Therefore, $x^{k}=x^{k+p}$ for all $p\geq 0$. Thus, the (GPPA) terminates after a finite number of steps. Without loss of generality, from now on, we assume that $f(x^k) >\ell^*$ for all $k$.

Recall that the (GPPA) starts from a point $x^0\in \mbox{\rm dom }g_1$ and generates two sequences $\{x^k\}$ and $\{y^k\}$ with $y^k\in \partial h(x^k)=\nabla h(x^k)$ and
$$y^{k-1} - \nabla g_2(x^{k-1}) - t(x^{k}- x^{k-1}) \in \partial^L g_1(x^{k}).$$
Thus, from Proposition \ref{sumrule} we have
\begin{align*}
\left(y^{k-1} - \nabla g_2(x^{k-1}) - t(x^{k}- x^{k-1})\right) + \nabla g_2(x^k) -y^k&\in \partial^L g_1(x^k) + \nabla g_2(x^k) - \nabla h(x^k) = \partial^Lf(x^k).
\end{align*}
Using the Lipschitz continuity of $\nabla g_2$ and $\nabla h$, we have
\begin{align*}
&\left \|y^{k-1} - \nabla g_2(x^{k-1})- t(x^{k}- x^{k-1}) + \nabla g_2(x^k) -y^k \right \| = \\
&= \left \| \left(\nabla h(x^{k-1}) - \nabla h(x^k) \right) + \left(\nabla g_2(x^k) - \nabla g_2(x^{k-1}) \right)- t(x^{k}- x^{k-1}) \right \|\\
& \leq \left(L(h)+ L + t \right)\|x^{k-1} -x^k\|  \leq M \|x^{k-1} -x^k\|,
\end{align*}
where $M:=L(h)+ L + t$. Therefore,
\begin{equation}
\mbox{dist}\left(0; \partial^Lf(x^k)\right) \leq M \|x^{k-1}-x^{k}\|.
\label{dist1}
\end{equation}
According to the assumption that $f$ has the strong Kurdyka - \L ojasiewicz property at $x^*$, there exist $\nu>0$, a neighborhood $\mathcal{V}$ of $x^*$, and  a continuous concave function $\varphi: \; [0, \nu[ \to [0, +\infty[$ so that for all $x\in \mathcal{V}$ satisfying $\ell^*<f(x) < \ell^*+\nu$, we have
\begin{equation}
\varphi'\left(f(x)-\ell^*\right)\mbox{dist}\left(0; \partial^L f(x)\right)\geq 1.
\label{KL1}
\end{equation}
Let $\delta>0$ small enough such that $\B(x^*; \delta) \subset \mathcal{V}$. Using the facts that $\lim \limits_{\ell\to +\infty} x^{k_\ell} =x^*$, $\lim \limits_{k\to +\infty}\|x^{k+1}-x^k\|=0$, $\lim \limits_{k\to +\infty} f(x^k) = \ell^*$, and $f(x^k)>\ell^*$ for all $k$, we can find a natural number $N$ large enough satisfying
\begin{equation}
x^N\in \B(x^*; \delta), \;\; \ell^*< f(x^N) < \ell^* +\nu,
\label{eq1}
\end{equation}
and
\begin{equation}
\|x^N-x^*\| + \frac{\|x^N-x^{N-1}\| }{4}+ \gamma\varphi\left(f(x^N) -\ell^*\right)<\frac{3\delta}{4},
\label{eq2}
\end{equation}
where $\gamma=\frac{2M}{t-L}>0$. We will show that for all $k\geq N$, $x^k \in \B(x^*; \delta)$. To this end, we first show that whenever $x^k\in \B(x^*; \delta)$ and $\ell^* < f(x^k) < \ell^*+\nu$ for some $k$, we have
\begin{align}\label{eq7}
\|x^{k} - x^{k+1}\| &\leq \frac{\|x^{k-1} - x^{k} \| }{4}  + \gamma\left[\varphi\left(f(x^k) -\ell^*\right) - \varphi\left(f(x^{k+1}) -\ell^*\right)\right].\end{align}
Indeed, by \eqref{dist1}, the concavity of $\varphi$, \eqref{KL1}, and \eqref{descentf}, we have
\begin{align*}
M\|x^{k-1}-x^{k}\|&\left[\varphi\left(f(x^k) -\ell^*\right) - \varphi\left(f(x^{k+1}) -\ell^*\right)\right] &\\
&\geq \mbox{dist}\left(0; \partial^Lf(x^k)\right) \left[\varphi\left(f(x^k) -\ell^*\right) - \varphi\left(f(x^{k+1}) -\ell^*\right)\right]\\
&\geq \mbox{dist}\left(0; \partial^Lf(x^k)\right) \varphi'\left(f(x^k)-\ell^*\right)\left[f(x^k) - f(x^{k+1})\right]\\
&\geq \frac{t-L}{2}\|x^{k} -x^{k+1}\|^2.
\end{align*}
It follows that
\begin{align}\label{eq3}
\varphi\left(f(x^k) -\ell^*\right) - \varphi\left(f(x^{k+1}) -\ell^*\right) & \geq  \frac{1}{\gamma}\frac{\|x^{k} -x^{k+1}\|^2}{\|x^{k-1}-x^k\|} \notag \\
&\geq \frac{1}{\gamma} \left[\|x^{k} - x^{k+1} \|  - \frac{\|x^{k-1} - x^{k} \| }{4}\right],
\end{align}
where the last inequality holds since $\frac{a^2}{b} \geq a - \frac{b}{4}$ for any positive real numbers $a$ and $b$. This implies \eqref{eq7}.

We next show that $x^k \in \B(x^*; \delta)$ for all $k\geq N$ by induction. The claim is true for $k=N$ by the construction above. Now suppose the assertion holds for $k=N, \ldots, N+k-1$ for some $k\geq 1$, i.e., $x^N, \ldots, x^{N+k-1}\in \B(x^*; \delta)$. Since $f(x^k)$ is a non-increasing sequence that converges to $\ell^*$, our choice of $N$ implies that $\ell^*< f(x^k) <\ell^*+\nu$ for all $k\geq N$. In particular, \eqref{eq7} can be applied for all $k=N, \ldots, N+k-1$. Using the estimation \eqref{eq7} for $k=N, \ldots, N+k-1$, we have
\begin{align*}
\|x^{N} - x^{N+1}\| &\leq \frac{\|x^{N-1} - x^{N} \| }{4}  + \gamma\left[\varphi\left(f(x^N) -\ell^*\right) - \varphi\left(f(x^{N+1}) -\ell^*\right)\right],\\
\|x^{N+1} - x^{N+2}\| &\leq \frac{\|x^{N} - x^{N+1} \| }{4}  + \gamma\left[\varphi\left(f(x^{N+1}) -\ell^*\right) - \varphi\left(f(x^{N+2}) -\ell^*\right)\right],\\
&\ldots\\
\|x^{N+k-1} - x^{N+k}\| &\leq \frac{\|x^{N+k-2} - x^{N+k-1} \| }{4}  + \gamma\left[\varphi\left(f(x^{N+k-1}) -\ell^*\right) - \varphi\left(f(x^{N+k}) -\ell^*\right)\right].
\end{align*}
Therefore,
\begin{align*}
\sum_{j=1}^k\|x^{N+j}-x^{N+j-1}\|  \leq \frac{1}{4}\sum_{j=1}^k\|x^{N+j}-x^{N+j-1}\| &+ \frac{\|x^{N-1} - x^{N} \| }{4}  - \frac{\|x^{N+k-1} - x^{N+k} \| }{4} \\
&+ \gamma\left[\varphi\left(f(x^N) -\ell^*\right) - \varphi\left(f(x^{N+k}) -\ell^*\right)\right].
\end{align*}
Making use of the non-negativity of $\varphi$, we get
\begin{equation}
\sum_{j=1}^k\|x^{N+j}-x^{N+j-1}\| \leq  \frac{4}{3}\left[\frac{\|x^{N-1} - x^{N} \| }{4} + \gamma\varphi\left(f(x^N) -\ell^*\right)\right].
\label{eq4}
\end{equation}
It follows that
\begin{align*}
\|x^{N+k} -x^*\| &\leq \|x^N-x^*\| +\sum_{j=1}^k\|x^{N+j}-x^{N+j-1}\| \\
&\leq \frac{4}{3}\left[\|x^N-x^*\| +  \frac{\|x^{N-1} - x^{N} \| }{4} + \gamma\varphi\left(f(x^N) -\ell^*\right)        \right] <\delta.
\end{align*}
Thus, $x^k \in \B(x^*; \delta)$ for all $k\geq N$. Since $x^k\in \B(x^*, \delta)$ and $\ell^*< f(x^k) <\ell^*+\nu$ for all $k\geq N$, it follows from~\eqref{eq4} by letting $k\to +\infty$ that $\sum_{k=1}^{+\infty}\|x^{k+1} -x^k\| <+\infty$. Therefore, $\{x^k\}$ is a Cauchy sequence and hence it is a convergent sequence. \qed

Below is another theorem which gives sufficient conditions that guarantee the convergence of the sequence $\{x_k\}$ generated by (GPPA). In contrast to Theorem \ref{main}, we require the differentiability with Lipschitz gradient of the function $g_1+g_2$ instead of $h$ along with the strong Kurdyka - \L ojasiewicz property of $f$. In this case, without loss of generality, we can assume that $g_1(x)=0$. In the next result, for convenience, we put $g_2(x) = g(x)$.

\begin{Theorem}\label{main2}
 Consider the difference of functions $f=g-h$  with $\inf_{x\in \R^n} f(x) >-\infty$. Suppose that $g$ is differentiable and $\nabla g$ is $L$ - Lipschitz continuous, $f$ has the strong Kurdyka - \L ojasiewicz property at any point $x\in \mbox{\rm dom} f$, and $h$ is a finite convex function.  If $C^*\neq \emptyset$, then the {\rm(GPPA)} sequence $\{x^k\}$ converges to a critical point of $f$.
\end{Theorem}
{\it Proof}  The proof is very similar to that of Theorem \ref{main}, except a few adjustments. Note that $f$ is locally Lipschitz continuous under the assumptions made since $g$ is a $C^1$ function and $h$ is a finite convex function. By \eqref{findyk}, we have
%$$y^{k-1} - \nabla g_2(x^{k-1}) - \left(\frac{1}{\lambda}+ L\right)(x^{k}- x^{k-1}) \in \partial g_1(x^{k}),$$
%and
%$$y^{k} - \nabla g_2(x^{k}) - \left(\frac{1}{\lambda}+ L\right)(x^{k+1}- x^{k}) \in \partial g_1(x^{k+1}).$$
$$y^{k-1} - t(x^k -x^{k-1}) = \nabla g(x^k ) \mbox{ and } y^{k} - t(x^{k+1} -x^{k}) = \nabla g(x^{k+1} ).$$
This implies,
$$y^{k} - \left(y^{k-1} - t(x^k -x^{k-1})   \right) = \nabla g(x^{k+1} ) -  \nabla g(x^k ) + t(x^{k+1} -x^{k}) .$$
Making use of the Lipschitz continuity of $\nabla g$ yields
\begin{align*}
\left\|y^{k} - \left(y^{k-1} - t(x^k -x^{k-1})   \right)\right\| &= \left \|\nabla g(x^{k+1} ) -  \nabla g(x^k ) + t(x^{k+1} -x^{k})\right\|\\
&\leq \left(L+t\right)\|x^{k} -x^{k+1}\|.
\end{align*}
On the other hand,
$$y^{k} - \left(y^{k-1} - t(x^k -x^{k-1})   \right) \in \partial h(x^k) -\nabla g(x^k) =\partial^F(-f)(x^k) \subset \partial^C(-f)(x^k).$$
Since $\partial^C(-f)(x^k) = -\partial^Cf(x^k)$, we have
$$\mbox{dist}\left(0; \partial^Cf(x^k)\right) =\mbox{dist}\left(0; \partial^C (-f)(x^k)\right) \leq \left(L+t\right)\|x^{k} -x^{k+1}\|.$$
Choose $N$ as in \eqref{eq1} and~\eqref{eq2} with $\gamma=\frac{2L+2t}{t-L}$ instead of $\frac{2M}{t-L}$ as before.  For all $k$ large enough such that $x^k\in \B(x^*; r)$ and $\ell^* < f(x^k) < \ell^*+\nu$, we have
\begin{align*}
\left(L+t\right)\|x^{k}-x^{k+1}\|&\left[\varphi\left(f(x^k) -\ell^*\right) - \varphi\left(f(x^{k+1}) -\ell^*\right)\right] &\\
&\geq \mbox{dist}\left(0; \partial^Cf(x^k)\right) \left[\varphi\left(f(x^k) -\ell^*\right) - \varphi\left(f(x^{k+1}) -\ell^*\right)\right]\\
&\geq \mbox{dist}\left(0; \partial^Cf(x^k)\right) \varphi'\left(f(x^k)-\ell^*\right)\left[f(x^k) - f(x^{k+1})\right]\\
&\geq \frac{t-L}{2}\|x^{k} -x^{k+1}\|^2.
\end{align*}
It follows that
\begin{equation}
\|x^{k} -x^{k+1}\| \leq \gamma \left[ \varphi\left(f(x^k) -\ell^*\right) - \varphi\left(f(x^{k+1}) -\ell^*\right)\right].
\label{est1}
\end{equation}
From this, the induction to prove that $x^k \in \B(x^*; r)$ for all $k\geq N$ can be carried out similarly to the proof of Theorem~\ref{main}. Indeed, suppose the assertion holds for $k=N, \ldots, N+k-1$ for some $k\geq 1$, i.e., $x^N, \ldots, x^{N+k-1}\in \B(x^*; r)$. Observe that
\begin{align*}
\|x^{N+k} -x^*\| &\leq \|x^N-x^*\| +\sum_{j=1}^k\|x^{N+j-1}-x^{N+j}\| \\
&\leq \|x^N-x^*\|  + \gamma \sum_{j=1}^k\left[ \varphi\left(f(x^{N+j-1}) -\ell^*\right) - \varphi\left(f(x^{N+j}) -\ell^*\right)\right]\\
&\leq \|x^N-x^*\|  + \gamma \varphi\left(f(x^N) -\ell^*\right) <r.
\end{align*}
Thus, $x^k \in \B(x^*; r)$ for all $k\geq N$. Since $x^k\in \B(x^*, r)$ and $\ell^*< f(x^k) <\ell^*+\nu$ for all $k\geq N$, we can sum~\eqref{est1} from $k=N$ to some $N_1$ greater than $N$ and take the limit as $N_1\to +\infty$, showing that $\sum_{k=1}^\infty\|x^{k+1} -x^k\| <+\infty$. This completes the proof. \qed

In the proposition below, we give sufficient conditions for the set of cluster points $C^*$ of the (GPPA) sequence $\{x_k\}$ to be nonempty.

\begin{Proposition} Consider the function $f = g - h$, where $g=g_1+g_2$ in \eqref{generalDC}. Let $\{x^k\}$ be sequence generated by the {\rm(GPPA)} for solving~\eqref{d.c.}. The set of critical points $C^*$ of $\{x_k\}$ is nonempty if one of the following conditions is satisfied:\\[1ex]
{\rm\bf (i)} For any $\alpha$, the lower level set $L_{\leq \alpha}:=\{x\in \R^n:\; f(x)\leq \alpha\}$ is bounded.\\
{\rm\bf (ii)} $\liminf\limits_{\|x\|\to +\infty} h(x)=+\infty$ and $\liminf\limits_{\|x\|\to +\infty} \frac{g(x)}{h(x)}>1$.
\end{Proposition}
{\it Proof} The conclusion under {\bf(i)} follows directly form the facts that $f(x^k)\leq f(x^0)$ for all $k$ and $L_{\leq f(x^0)}$ is bounded. Now assume that {\bf(ii)} is satisfied. Then there exist $M>1$ and $R>0$ such that $g(x) \geq Mh(x)$ for all $x$ satisfying $\|x\|\geq R$. It follows that
$$\liminf_{\|x\|\to +\infty}f(x) = \liminf_{\|x\|\to +\infty} \left[g(x) - h(x)\right] \geq (M-1)\liminf_{\|x\|\to +\infty} h(x) =+\infty.$$
Thus, $f$ is coercive. Combining this with the descent property of the sequence $\{f(x^k)\}$, we can conclude that $\{x^k\}$ is bounded.\qed

It is known from  \cite[Corollary 16]{BDLS07} and \cite[Section 4.3]{ABRS10} that a proper lower semicontinuous semi-algebraic function $f$ on $\R^n$ always satisfies the Kurdyka - \L ojasiewicz property at all points in $\mbox{dom }\partial f$ with
$\varphi(s) = cs^{1-\theta}$ for some $\theta \in [0, 1[$ and $c>0$. We now derive convergence rates
of the (GPPA) sequence by examining the range of the exponent.

\begin{Theorem} Consider the settings of Theorems \ref{main} and \ref{main2}. Suppose further that $f$ is a proper closed semi-algebraic function so that the function $\varphi$ in the  Kurdyka - \L ojasiewicz property has
 the form $\varphi(s) = cs^{1-\theta}$ for some $\theta \in [0, 1[$ and $c>0$.  Then we have the following conclusions.\\[1ex]
{\rm\bf (i)} If $\theta=0$, then the sequence $\{x^k\}$ converges in a finite number of steps.\\
{\rm\bf (ii)} If $0 < \theta \leq \frac{1}{2}$, then there exist $\mu>0$ and $q \in (0, 1)$ satisfying
$$\|x^k -x^*\|\leq \mu q^k.$$
{\rm\bf (iii)} If $\frac{1}{2}<\theta < 1$, then there exists $\mu>0$ such that
$$\|x^k - x^*\| \leq \mu k^{\frac{1-\theta}{1-2\theta}}.$$
\end{Theorem}
{\it Proof} For each $k \geq 1$, set $\Delta_k=\sum_{p=k}^{+\infty} \|x^{p+1}-x^p\|$ and set $\ell_k= f(x^k) -\ell^*$. It is obvious from the triangle inequality that $\|x^k -x^*\| \leq \Delta_k$. From Kurdyka - \L ojasiewicz property with the special form of $\varphi$, we have
\begin{equation}
c(1-\theta)\ell_k^{-\theta}\mbox{dist}\left(0; \partial^L f(x^k)\right) \geq 1.
\label{lower_est}
\end{equation}
From the proof of Theorem \ref{main2}, if $\nabla g$ is $L$ - Lipschitz continuous, then
$$\mbox{dist}\left(0; \partial^L f(x^k)\right) \leq \left(L+t\right) \|x^{k+1} - x^k\|,$$
for all sufficiently large $k$. Combining this with \eqref{est1} yields
$$\Delta_k \leq \gamma \varphi(\ell_k) \leq \gamma \varphi(\ell_{k-1}) = \gamma c\ell_{k-1}^{1-\theta}\leq \gamma c \left[(L+t)c(1-\theta)\right]^{\frac{1-\theta}{\theta}}\|x^k-x^{k-1}\|^{\frac{1-\theta}{\theta}},$$
where $\gamma=\frac{2L+2t}{t-L}$.

In the case of  Theorem \ref{main} where $\nabla h$ is $L(h)$ - Lipschitz continuous, we have
$$\mbox{dist}\left(0; \partial^L f(x^k)\right) \leq M \|x^{k} - x^{k-1}\|,$$
for all sufficiently large $k$, where $M=L(h)+L+t$. It follows from~\eqref{eq4} that
$$\Delta_k\leq \frac{4}{3}\left[\frac{\|x^k-x^{k-1} \|}{4} +\gamma\varphi(\ell_k)  \right] \leq \frac{\|x^k-x^{k-1} \|}{3} + \frac{4\gamma}{3}\left[Mc(1-\theta)\right]^{\frac{1-\theta}{\theta}}\|x^k-x^{k-1}\|^{\frac{1-\theta}{\theta}},$$
where $\gamma=\frac{2L+2t}{t-L}$. Thus, in both cases it always holds that
$$\Delta_k\leq C_1 (\Delta_{k-1} -\Delta_k) +C_2\left(\Delta_{k-1} -\Delta_k\right)^{\frac{1-\theta}{\theta}},$$
for some $C_1, C_2>0.$ The result now follows from the proof of \cite[Theorem 2]{AB09}. \qed

\section{Examples}
\setcounter{equation}{0}

\noindent {\bf Trust-Region SubProblem. } Consider the \emph{trust-region subproblem}
\begin{equation}\label{TRB}
\min\bigg\{\phi(x)=\frac{1}{2}x^\top Ax + b^\top x:\; \|x\|^2\leq r^2\bigg\},
\end{equation}
where $A$ is an $n\times n$ real symmetric matrix and $b\in \R^n$ is given. Since $A$ is not  required to be positive-semidefinite, \eqref{TRB} is  a nonconvex optimization problem. Let  $E=\{x\in \R^n: \; \|x\|\leq r\}$ and define the function
\begin{equation*}
f(x)=\phi(x)+\delta(x; E), \; x\in \R^n.
\end{equation*}
The trust-region subproblem \eqref{TRB} can be solved by the ${\rm(DCA)}$ with the following DC decomposition $f=g - h$ with
\begin{equation*}
g(x)=\frac{1}{2}\rho\|x\|^2+b^\top x+\delta(x; E)\; \mbox{\rm and }h(x)=\frac{1}{2}x^\top(\rho I-A)x,
\end{equation*}
where $\rho$ is a positive number such that $\rho I-A$ is positive-semidefinite; see \cite{TA2}. The convergence analysis of the (DCA) sequence for solving \eqref{TRB} was proved in \cite{TuanYen}.

Define
 $$g_2(x) = \frac{1}{2}\rho\|x\|^2+b^\top x \mbox{ and } g_1(x)=\delta(x; E).$$
In this case, $g_2$ and $h$ have Lipschitz gradient with Lipschitz constants $L=\rho$ and $L(h) = \lambda_{\max}(\rho I-A)$, respectively. Applying the (GPPA) for \eqref{TRB}, we have $y^k = \nabla h(x^k)=(\rho I-A)x^k$ and
$$y^k - \nabla g_2(x^{k+1})  - t \left(x^{k+1} -x^k\right) \in \partial g_1(x^{k+1}).$$
This implies
$$y^k +tx^k - b \in (t+\rho)x^{k+1} +N(x^{k+1}; E).$$
Thus,
$$x^{k+1}= P_E\left(  \frac{1}{t+\rho}\left( (t+\rho)x^k - Ax^k- b\right)\right).$$

\begin{Proposition}
Consider the trust-region subproblem \eqref{TRB}. Then $C^*\neq \emptyset$ and the {\rm (GPPA)} sequence $\{x^k\}$ converges to a critical point of $f=g_1+g_2-h$.
\end{Proposition}
{\it Proof} We only need to verify that all assumptions of Theorem \ref{main} are satisfied in this particular case. Note that $f(x)=\phi(x)+\delta(x;E)$. Obviously, $\inf_{x\in \R^n}f(x)>-\infty$ and $C^*\neq \emptyset$. Let us show that $f$ is a semi-algebraic function. Note that
$$E=\{x\in \R^n:\; p(x)\leq r^2\},$$
where $p$ is the polynomial $p(x)=\sum_{i=1}^nx_i^2$. Thus, $E$ is a semialgebraic set, which implies that its associated indicator function is a semi-algebraic function; see, e.g., \cite{ABRS10}.

It is also straightforward that $\phi$ is also a semi-algebraic function since its graph
$$\gph \phi=\{(x,y)\in \R^n\times \R:\; x^\top Ax+b^\top x-y=0\}$$
is a semi-algebraic set. It follows that $f$ is a semialgebraic function as it is the sum of two semi-algebraic functions; see, e.g., \cite{ABRS10}. Therefore, $f$ satisfies the Kurdyka - \L ojasiewicz property.  Obviously, $h$ has Lipschitz continuous gradient. We have shown that all assumptions of Theorem \ref{main} are satisfied and the conclusion follows from Theorem \ref{main}. \qed

\noindent {\bf Nonconvex Feasibility Problems.} In this part, we show how the (GPPA) can be applied to solve nonconvex feasibility problems. Let $A$ and $B$ be two nonempty closed sets in $\R^n$. It is implicitly assumed that $A$ and $B$ are simple enough so that the projection onto each set is easy to compute. The feasibility problem asks for a point in $A\cap B$. It is clear that $A\cap B \neq \emptyset$
if and only if the following optimization problem has the
zero optimal value:
\begin{equation}
\min \bigg\{\frac{1}{2}d^2_B(x):\; x\in A\bigg\}.
\label{feas_prob}
\end{equation}
This problem is of the type  \eqref{generalDC} with the objective function $f(x) = g_1(x) + g_2(x) -h(x)$, where
$$g_1(x) = \delta(x; A), \; g_2(x) = \frac{1}{2}\|x\|^2, h(x)=\frac{1}{2}\left(\|x\|^2 -d^2_B(x)\right).$$
Obviously, the function $g_2$ is differentiable with $L-$Lipschitz gradient where $L=1$. We have
\begin{align*}
h(x)& = \dfrac{1}{2}\|x\|^2 - \dfrac{1}{2}\inf\{\|x\|^2 + \|y\|^2 -2\la x,y\ra: \; y\in B\}\\
&=\sup\{ \la x,y\ra - \frac{\|y\|^2}{2}: \; y\in B\}\\
&=\sup\{f_y(x): \; y\in B\},
\end{align*}
where $f_y(x)=\la x,y\ra - \frac{\|y\|^2}{2}$. Therefore, $h$ is a pointwise supremum of  a collection of affine functions so it is a convex function. Denote $S(\ox) = \{y\in B:\; f_y(\ox)=h(\ox)\}$. We have
$$S(\ox)=\{y\in B:\; \|\ox-y\|^2 =d^2_B(\ox)\}=P_B(\ox).$$
Since $B$ is a nonempty closed subset of $\R^n$, the set $S(\ox)=P_B(\ox)$ is nonempty and compact for any $\ox\in \R^n$. By  \cite[Theorem 3, p. 201]{IT79}, we have
$$\partial h(\ox) = \overline{\mbox{co}} \left(\bigcup_{y\in S(\ox)}\partial f_y(\ox)\right) = \overline{\mbox{co}} \left(\bigcup_{y\in S(\ox)} \{y\}\right)= \mbox{co }P_B(\ox).$$
Making use of Proposition 3.1, we now can state the necessary condition for a local minimum of \eqref{feas_prob}.
\begin{Proposition}If  $\ox\in A$ is a local optimal solution of \eqref{feas_prob}, then
\begin{align}
P_B(\ox) \subset \ox + N^L(\ox; A),
\label{optcond}
\end{align}
where  $N^L(\ox; A)$ is the limiting normal cone to $A$ at $\ox$ defined by $N^L(\ox; A)=\partial^L\delta(\ox;A)$.
\end{Proposition}
Note that the optimality condition \eqref{optcond} is not sufficient to ensure that $\ox$ is a local minimizer of \eqref{feas_prob} as shown in the next example.
\begin{Example}{\rm
Consider the following subsets of $\R^2$:
$$A=\big\{(x_1, x_2) \in \R^2:\; x_2 \geq 1\big\}\; \mbox{\rm \; and } B=\big \{(x_1, x_2)\in \R^2:\; x_2\leq \alpha x_1^2\big\},$$
where $\alpha <\frac{1}{2}$. Put  $\ox=(0, 1) \in A$. Since $\alpha<\frac{1}{2}$, the system
$$\begin{cases}
x_1^2 +(x_2-1)^2 &\leq 1,\\
x_2 - \alpha x_1^2&\leq 0,
\end{cases}
$$
has a unique solution $(x_1, x_2)=(0, 0)$.
This implies $P_B(\ox)=\{(0,0)\}$ and $d_B(\ox)=1$. Obviously, $\ox$ satisfies condition~\eqref{optcond} since
$$P_B(\ox) =\{(0, 0)\} \subset \{(0, \gamma):\; \gamma \leq 1\} =\ox +N(\ox; A).$$
However, for any neighborhood $U$ of $\ox$, there always exists $\epsilon>0$ small enough such that $x_\epsilon=(\epsilon, 1) \in U$ and
$$d_B(x_\epsilon)  \leq 1-\alpha \epsilon^2 <1.$$
Thus, $\oz$ cannot be a local minimizer of \eqref{feas_prob}.
}
\end{Example}
Based on the (GPPA), we now propose the following simple algorithm for solving \eqref{feas_prob}. For a given initial point $x^0\in A$, the (GPPA) sequence $\{x^k\}$ with the starting point $x^0$ is defined by
\begin{equation}
x^{k+1}\in P_A\left( (1-\frac{1}{t}) x^k + \frac{1}{t}y^k \right),\label{dca_seq}
\end{equation}
where $y^k$ is an element chosen in $\mbox{\rm co }P_B(x^k)$. Note that, this scheme is different from some other well-known methods such as the alternating projection algorithm or the averaged projection algorithm. Moreover, it cannot be obtained from the proximal forward - backward schemes in \cite{ABS11,BST}.

\begin{Theorem}
Let $A$ and $B$ are nonempty closed sets in $\R^n$ and let  $t>1$. Then the sequence $\{x^k\}\subset A$ satisfies the following: \\[1ex]
(i) For any $k\geq 1$,
$$d^2_B(x^k) -d^2_B(x^{k+1}) \geq 2(t-1)\|x^k-x^{k+1}\|^2.$$
(ii)$ \lim \limits_{k\to +\infty}\|x^k -x^{k+1}\| =0.$\\
(iii) If $\{x^k\}$ is bounded, then every cluster point is a critical point of $f=\delta(\cdot; A)+d^2_B(\cdot)$.
 \end{Theorem}

\begin{Proposition}
Let $A$ and $B$ are nonempty closed sets in $\R^n$ such that both of them are semi-algebraic sets and $B$ is convex. Suppose further that either $A$ or $B$ is bounded. Then the sequence $\{x^k\}$ generated by the {\rm (GPPA)} converges to a critical point of \eqref{feas_prob}.
\end{Proposition}
{\it Proof }As $A$ is a semi-algebraic set, the indicator function $\delta(\cdot; A)$ is a semi-algebraic function. On the other hand, $B$ is also semi-algebraic, so $x\mapsto \frac{1}{2}d^2_B(x)$ is also a semi-algebraic function; see \cite[Lemma 2.3]{ABS11}. Therefore, $f(x)=\delta(x; A) +\frac{1}{2}d^2_B(x)$ is a semi-algebraic function. If $B$ is closed and convex, it is well known that the function $x\mapsto d^2_B(x)$ is smooth with
1 - Lipschitz continuous gradient; see \cite[Corollary 12.30]{BausComb11}. The result now follows directly from Theorem \ref{main} since the boundedness of $\{x^k\}$ is ensured by the coercivity of $f$ under the assumption that either $A$ or $B$ is bounded. \qed

\section{Concluding Remarks}

Based on recent progress in using the Kurdyka - \L ojasiewicz property and variational analysis in analyzing nonsmooth optimization algorithms, we introduce and study convergence analysis of a proximal point algorithm for minimizing differences of functions. We are able to relax some convexity in the classical DC programming to deal with a more general class of problems. The results open up the possibility of understanding the convergence of the (DCA) and other algorithms for minimizing differences of convex functions used in numerous applications.

{\bf Acknowledgment. } The authors are very grateful to Prof. J\'{e}r\^ome Bolte for his helpful suggestions to improve the paper. The authors are also   thankful to Prof. Nguyen Dong Yen and Dr. Hoang Ngoc Tuan for useful discusions on the subject.  This work was completed while the first
author was visiting the Vietnam Institute for Advanced Study in Mathematics (VIASM). He would like to thank the VIASM for financial support and hospitality. The research of the second author was partially supported by the USA National Science Foundation under grant DMS-1411817 and
the Simons Foundation under grant \#208785.

\end{document}